\def \R {\hat{R}}

\def \be{\begin{equation}}
\def \eq{\end{equation}}
\def \ba{\begin{eqnarray}}
\def \ea{\end{eqnarray}}
\def \nn {\nonumber}
\def \pa {\partial}
\def \xh {\hat{x}}
\def \ph {\hat{\partial}}
\def \xih {\hat{\xi}}

\documentstyle[12pt]{article}

\begin{document}
\begin{titlepage}
\begin{center}
Jan 24, 1995
        \hfill  LBL-36746 \\
          \hfill    UCB-PTH-95/04  \\
\vskip .3in
{\large REALIZATION OF VECTOR FIELDS FOR QUANTUM GROUPS AS PSEUDODIFFERENTIAL
OPERATORS ON QUANTUM SPACES}
\footnote{This work was supported in part by the Director, Office of
Energy Research, Office of High Energy and Nuclear Physics, Division of
High Energy Physics of the U.S. Department of Energy under Contract
DE-AC03-76SF00098 and in part by the National Science Foundation under
grant PHY-90-21139.}
\vskip .3in
{CHONG-SUN CHU\footnote{email: cschu@physics.berkeley.edu}
 and BRUNO ZUMINO
}
\vskip .2in
{\em Department of Physics \\
  and\\
   Theoretical Physics Group\\
    Lawrence Berkeley Laboratory\\
      University of California\\
    Berkeley, CA 94720\\
   USA}
\end{center}
\vskip .3in
\begin{abstract}
The vector fields of the quantum Lie algebra are described for the quantum
groups $GL_q(N), SL_q(N)$
and $SO_q(N)$
as pseudodifferential operators on the linear quantum spaces covariant under
the corresponding quantum group. Their
expressions are simple and compact. It is pointed out that these vector fields
satisfy certain characteristic polynomial identities. The real forms $SU_q(N)$
and $SO_q(N,R)$ are discussed in detail.
\end{abstract}
\end{titlepage}
\renewcommand{\thepage}{\roman{page}}
\setcounter{page}{2}
\mbox{ }

\vskip 1in

\begin{center}
{\bf Disclaimer}
\end{center}

\vskip .2in
\begin{scriptsize}
\begin{quotation}
This document was prepared as an account of work sponsored by the United
States Government. While this document is believed to contain correct
 information, neither the United States Government nor any agency
thereof, nor The Regents of the University of California, nor any of their
employees, makes any warranty, express or implied, or assumes any legal
liability or responsibility for the accuracy, completeness, or usefulness
of any information, apparatus, product, or process disclosed, or represents
that its use would not infringe privately owned rights.  Reference herein
to any specific commercial products process, or service by its trade name,
trademark, manufacturer, or otherwise, does not necessarily constitute or
imply its endorsement, recommendation, or favoring by the United States
Government or any agency thereof, or The Regents of the University of
California.  The views and opinions of authors expressed herein do not
necessarily state or reflect those of the United States Government or any
agency thereof of The Regents of the University of California and shall
not be used for advertising or product endorsement purposes.
\end{quotation}
\end{scriptsize}

\vskip 2in

\begin{center}
\begin{small}
{\it Lawrence Berkeley Laboratory is an equal opportunity employer.}
\end{small}
\end{center}
\newpage

\renewcommand{\thepage}{\arabic{page}}
\setcounter{page}{1}
\section {Vector Fields for Quantum Groups}
A quantum group can be described\cite{RTF} in terms of matrices $A$ with
noncommuting
elements satisfying the equation
\be  \label{RAA} \R_{12} A_1 A_2 =A_1 A_2 \R_{12}, \eq
with the $\R$ matrix appropriate to the particular quantum group.
The matrix elements generate the algebra of functions on the group.
Here we have used a well known standard notation: for instance, Eq.(\ref{RAA})
written explicitly, takes the form
\be \R^{ij}_{kl} A^k_m  A^l_n =A^i_k A^j_l \R^{kl}_{mn}. \eq
The vector fields on the quantum group can be described \cite{Z1,Z2,SWZ1,SWZ2}
by the matrix elements of a matrix $Y$ satisfying the commutation relation
\be \label{YY} \R_{12} Y_2 \R_{12} Y_2 =Y_2 \R_{12} Y_2 \R_{12}, \eq
which corresponds to the Lie algebra relations in the classical case. The
action of the vector fields on the group is then given by the commutation
relation
\be Y_1 A_2 =A_2 \R_{12} Y_2 \R_{12}. \eq

The quantum group matrices can coact on a quantum space, for instance by right
multiplication. A point of coordinates $x_{0i}$ will be transformed into
$x_i=x_{0j} A^j_i$ or, more compactly,
\be x=x_0 A. \eq
Keeping the original point $x_0$ fixed, the action of a vector field on the
quantum group induces an action on the quantum space
\be \label{Yx} Y_1 x_2 =x_2 \R_{12} Y_2 \R_{12}, \eq
i.e. \be Y^i_j x_k = x_m \R^{im}_{ln} Y^n_r \R^{lr}_{jk}. \eq

We shall consider the case when a differential calculus covariant with respect
to the coaction of the quantum group exists on the quantum space.
In this case it is natural to ask whether it is possible to realize the vector
fields $Y$
as pseudodifferential operators satisfying Eqs.(\ref{YY}) and (\ref{Yx}). We
shall show that
this can be done for the quantum groups $GL_q(N), SL_q(N)$ and $SO_q(N)$. Their
real forms are also
considered.

\section{$GL_q(N), SL_q(N)$ and $SU_q(N)$}
The calculus for the quantum plane covariant under $GL_q(N)$ is well
known\cite{WZ}.
The coordinates $x_i$ in the plane satisfy the commutation relations
\be \label{xx} x_1 x_2 =q^{-1} x_1 x_2 \R_{12} \eq
and the derivatives $\pa^i$ satisfy
\be \label{px} \pa^i x_j = \delta^i_j +q \R^{ik}_{jl} x_k \pa^l \eq
and
\be \label{pp} \pa_2 \pa_1 =q^{-1} \R_{12} \pa_2 \pa_1. \eq
All indices run from 1 to $N$ and $\R$ is the $GL_q(N)$ matrix, which satisfies
the characteristic equation
\be \label{char} \R^2 =1+\lambda \R, \quad \lambda=q-q^{-1}. \eq
Using Eq.(\ref{char}) and the above commutation relations, it is easy to verify
that the differential operator
\be \label{Y-GL} Y^i_j =q^{-2} \delta^i_j +q^{-1}\lambda \pa^i x_j \eq
satisfies Eq.(\ref{Yx}), which we repeat here,
\be  Y_1 x_2 =x_2 \R_{12} Y_2 \R_{12}, \eq
as well as
\be \label{Yp} \pa_2 Y_1 =\R_{12} Y_2 \R_{12} \pa_2. \eq
Combining these two, one finds that the matrix $Y$ satisfies also
Eq.(\ref{YY}).

The quantum subgroup $SL_q(N)$ can be obtained from $GL_q(N)$ as
follows\cite{SWZ1}.
For the quantum matrices one uses the standard quantum determinant $det_qA$ and
defines a new matrix
\be \label{TA} T=(det_q A)^{-1/N} A \eq
having quantum determinant equal to one. For the vector fields, one defines an
appropriate determinant $Det Y$
and defines a new matrix of vector fields\cite{SWZ1,SWZ2}
\be Z=(DetY)^{-1/N} Y\eq
having determinant one.
The number of independent elements of the matrix $Z$ is $N^2-1$, as in the
classical case. For the particular
representation Eq.(\ref{Y-GL}) of the $Y$ matrix, it is possible to show that
\be Det Y =\mu, \eq
where $\mu$ is the rescaling operator in the plane
\be \label{mu} \mu=1+q \lambda x_i \pa^i, \eq
which satisfies
\be \mu x_i =q^2 x_i \mu, \quad \pa^i \mu =q^2 \mu \pa^i. \eq
Thus here
\be \label{ZY} Z=\mu^{-1/N} Y \eq
realizes the $SL_q(N)$ vector fields as pseudodifferential operators in the
quantum plane. Note that $\mu$ commutes with the elements of $Y$.

It is very easy to verify that the matrix given by Eq.(\ref{Y-GL}) satisfies
the identity
\be (Y-\mu)(Y-q^{-2})=0, \eq
where matrix multiplication is implied. This is a special example of polynomial
characteristic equations satisfied by quantum vector fields\cite{GZB}. In
general these
equations are of higher order but for the realization Eq.(\ref{Y-GL}) we see
that
the polynomial is quadratic in $Y$.
We intend to come back to a general treatment of these characteristic equations
in a forthcoming publication.

The invariant quantum trace of the k-th power of the matrix $Y$ is defined as
\be t_k=Tr D^{-1} Y^k, \eq
where $D$ is the diagonal matrix $(1,q^2,...,q^{2(N-1)})$. The $t_k$ commute
with
the matrix elements of $Y$. In general only the first $r$ $(k=1,2,...,r)$ are
independent, where $r$ is the rank of the
group\cite{RTF}, a fact which is related to the existence of the characteristic
polynomial equations for $Y$ mentioned above.
For $Y$ given by Eq.(\ref{Y-GL}) all $t_k$ are simply functions of $\mu$. For
instance,
\be t_1=[N]-1+\mu=q^{-2} t_0-q^{-2N} +\mu, \eq \nn
\be t_2=q^{-2} t_1-\mu q^{-2N} +\mu^2, \eq \nn
\be t_3=q^{-2} t_2 -\mu^2 q^{-2N} +\mu^3, \eq
etc., where
\be [N]=1+q^{-2}+q^{-4}+...+q^{-2(N-1)}. \eq

If $|q|=1$ the calculus given by Eqs.(\ref{xx}-\ref{pp}) for the quantum plane
can be given
a reality structure\cite{WZ,Sch} by requiring $x_i$ to be real
\be \overline{x_i} =x_i \eq
and by defining conjugate derivatives as
\be \overline{\pa^i} =-q^{2i'} \pa^i, \eq
where we have introduced the notation
\be \label{i'} i'=N+1-i, \quad i=1,2,...,N. \eq
Here we consider instead the case when $q$ is real and the complex conjugates
of $x_i$ and of $\pa^i$ are new independent variables. It will be convenient to
give them new names, i.e. we set
\be \overline{x_i}=\xh^i \eq
and
\be \overline{\pa^i}=-\ph_i. \eq
The commutation relations of these new variables can be obtained immediately
from Eqs.(\ref{xx}-\ref{pp}) by complex conjugation (remember that this is an
involution which inverts the order of factors in a product). Using the symmetry
property
\be \R^{ij}_{kl}=\R^{kl}_{ij}, \eq
we see that
\be \xh_2 \xh_1=q^{-1} \R_{12} \xh_2 \xh_1, \eq
\be \label{xhph}\xh^j \ph_i = -\delta^j_i +q \R^{jl}_{ik} \ph_l \xh^k \eq
and
\be \ph_1 \ph_2 =  q^{-1} \ph_1 \ph_2 \R_{12}. \eq
Eq.(\ref{xhph}) can be written in a form  more analogous to Eq.(\ref{px}) if
one introduces the matrix
\be \Psi^{ir}_{js} =(\R^{-1})^{ri}_{sj} q^{2(j-r)}
=(\R^{-1})^{ri}_{sj} q^{2(i-s)}, \eq
which satisfies
\be \R^{kj}_{li} \Psi^{ir}_{js} = \Psi^{kj}_{li} \R^{ir}_{js}
=\delta^k_s \delta^r_l \eq
and
\be \Psi^{ri}_{si} = \delta^r_s q^{-2(N-r)-1}, \eq
\be \Psi^{ir}_{is} = \delta^r_s q^{-2(r-1)-1}. \eq
It takes the form
\be \hat{\pa}_i \hat{x}^j= \delta^j_i q^{-2i'}+q^{-1} \Psi^{jl}_{ik} \hat{x}^k
\hat{\pa}_l, \eq
where $i'$ is given by Eq.(\ref{i'}).

To complete the algebra of the complex calculus, we must now give commutation
relations between the variables $x_i, \pa^i$ and their conjugate $\xh^i,
\ph_i$. A consistent set is given by
\be \xh^i x_j=q (\R^{-1})^{ik}_{jl} x_k \xh^l, \eq
\be \pa^i \xh^j = q (\R^{-1})^{ji}_{lk} \xh^k \pa^l, \eq
\be \ph_i x_j =q^{-1} \R^{kl}_{ij} x_k \ph_l \eq
and
\be \pa^i \ph_j =q^{-1} \R^{ik}_{jl} \ph_k \pa^l. \eq
Consistency can be checked by verifying that all these relations braid
correctly with each other.

Having the complex calculus we can now ask how the vector field realization of
Eq.(\ref{Y-GL}) acts on the conjugate variables. It is not hard to verify that
\be \label{Yxh} \xh_2 Y_1 = \R_{12}Y_2 \R^{-1}_{12} \xh_2 \eq
and
\be \label{Yph}  Y_1 \ph_2 =\ph_2 \R_{12} Y_2 \R^{-1}_{12}. \eq
On the other hand, by complex conjugation, Eqs.(\ref{Yx}),(\ref{Yp}),
(\ref{Yxh}) and (\ref{Yph}) give
\be \xh_2 Y_1^\dagger = \R_{12}Y_2^\dagger \R_{12} \xh_2, \eq
\be Y_1^\dagger \ph_2 =\ph_2 \R_{12} Y_2^\dagger \R_{12}, \eq
\be  Y_1^\dagger x_2 =x_2 \R^{-1}_{12} Y_2^\dagger \R_{12} \eq
and
\be \pa_2 Y_1^\dagger =\R^{-1}_{12} Y_2^\dagger \R_{12} \pa_2, \eq
where $Y^\dagger$ is the hermitian conjugate of the matrix $Y$
\be (Y^\dagger)^i_j =\overline{Y_i^j}
= q^{-2}\delta^i_j -q^{-1} \lambda \xh^i \ph_j,
\eq
which satisfies the equation conjugate of Eq.(\ref{YY})
\be \label{YhYh} \R_{12} Y_2^\dagger \R_{12} Y_2^\dagger =Y_2^\dagger \R_{12}
Y_2^\dagger \R_{12}, \eq
as well as the commutation relation with $Y$
\be \label{YYh} \R_{12} Y_2 \R^{-1}_{12} Y_2^\dagger =Y_2^\dagger \R_{12} Y_2
\R^{-1}_{12}. \eq

Until now, we have considered two $GL_q(N)$ groups complex conjugate of each
other,
i.e. a truly complex $GL_q(N)$\cite{DSWZ,CDSWZ,Pod}.
The quantum group can be restricted to $U_q(N)$ by imposing on its matrices the
unitarity condition
\be A^\dagger= A^{-1} \eq
and to $SU_q(N)$ by further normalizing the matrices as in Eq.(\ref{TA}) so
that they have quantum
determinant equal to one.

The vector fields of the $U_q(N)$ subgroup can be defined as the elements of
the Hermitian matrix
\be U=Y Y^\dagger. \eq
Indeed, it is very easy to check that $U$ commutes with the Hermitian length
\be {\cal L}=x_i \xh^i =x_i\overline{x_i} \eq
($Y$ and $Y^\dagger$ separately do not), i.e. the $U$ vector fields leave
${\cal L}$ invariant. $U$ is a perfectly
good matrix of vector fields and satisfies equations similar to Eq.(\ref{YY})
and Eq.(\ref{Yx})
\be \R_{12} U_2 \R_{12} U_2 = U_2 \R_{12} U_2 \R_{12}, \eq
\be U_1 x_2 =x_2 \R_{12} U_2 \R_{12} \eq
and
\be \xh_2 U_1 =\R_{12} U_2 \R_{12} \xh_2, \eq
as a consequence of equations for $Y$ and $Y^\dagger$ given above.
Notice that
\be \label{Uxp} q^2 U^i_j=q^{-2} \delta^i_j +q^{-1}\lambda \pa^i x_j -q^{-1}
\lambda \xh^i \ph_j -\lambda^2 \pa^i {\cal L} \ph_j \eq
which will be useful to us later.

Finally we observe that, if we want to reduce the  vector fields to the number
appropriate to $SU_q(N)$, we must normalize $U$,
i.e., take the matrix
\be \label{ZZU} Z Z^\dagger =U/(\mu \bar{\mu})^{1/N}. \eq
In addition to commuting with $Y^i_j$, the rescaling operator $\mu$ in
Eq.(\ref{mu}) commutes with
$\xh^i, \ph_i$ and therefore with $(Y^\dagger)^i_j$ and
\be \bar{\mu}=1-q\lambda \ph_i \xh^i. \eq
On the other hand $\bar{\mu}$ commutes with $(Y^\dagger)^i_j, x_i, \pa^i,
Y^i_j$ and satisfies
\be \bar{\mu} \xh^i=q^{-2} \xh^i \bar{\mu},\quad \ph_i \bar{\mu}=q^{-2}
\bar{\mu}\ph_i. \eq
Clearly $\mu \bar{\mu}$ commutes with ${\cal L}$, therefore so does
$ZZ^\dagger$. $Z$ and $Z^\dagger$
satisfy equations analogous to Eq.(\ref{YY}),(\ref{YhYh}),(\ref{YYh}). Using
this fact
one can show that
\be Det ZZ^\dagger =(DetZ)(DetZ^\dagger)=1. \eq
Notice that the vector field matrix $Z Z^\dagger$ is Hermitian, which is the
natural
reality condition for
$SU_q(N)$.

\section{$SO_q(N)$ and $SO_q(N,R)$}
We shall call $T$ the quantum matrices of $SO_q(N)$, instead of $A$. In
addition to
\be \R_{12} T_1 T_2 =T_1 T_2 \R_{12}, \eq
they satisfies the orthogonality relations\cite{RTF}
\be T^t g T=g, \quad T g^{-1} T^t = g^{-1}, \eq
where the numerical quantum metric matrices $g=g_{ij}$ and $g^{-1}=g^{ij}$ can
be
chosen to be equal $g_{ij}=g^{ij}$. The $SO_q(N)$ $\R$ matrix satisfies also
orthogonality
conditions
\be (\R^{-1})^{ij}_{kl}=g^{im} \R^{jn}_{mk} g_{nl}=g_{km} \R^{mi}_{ln} g^{nj},
\eq
as well as the usual symmetry relations
\be \R^{ij}_{kl} =\R^{kl}_{ij}. \eq

The $SO_q(N)$ vector field matrix, which we shall call $Z$, satisfies
\be \label{ZZ} \R_{12} Z_2 \R_{12} Z_2 = Z_2 \R_{12} Z_2 \R_{12}, \eq
\be \label{ZT} Z_1 T_2 =T_2 \R_{12} Z_2 \R_{12}, \eq
as well as an orthogonality constraint in one of the two equivalent
forms\cite{Z2,SWZ2}
\be \label{gZZ} g_{ij} (Z_2 \R_{12} Z_2)^{ij}_{kl} =q^{1-N} g_{kl}, \eq
\be \label{ZZg} (Z_2 \R_{12} Z_2)^{ij}_{kl} g^{kl} =q^{1-N} g^{ij}. \eq
Eq.(\ref{gZZ}) or (\ref{ZZg}) reduces the number of independent vector fields
from
$N^2$ to $N(N-1)/2$ as in the classical case.

The projector decomposition of the $\R$ matrix for $SO_q(N)$ is
\be \R =q P^+-q^{-1} P^{-} +q^{1-N}P^0. \eq
Here $P^+$ is the traceless part of the symmetriser,
$P^-$ is the antisymmetriser and $P^0$ is the trace operator.
It is related to the metric by
\be (P^0)^{ij}_{kl}=\nu g^{ij} g_{kl}, \quad \nu =\frac{\lambda}{(q^N
-1)(q^{1-N}+q^{-1})}. \eq
The coordinates $x_i$ of the quantum Euclidean space satisfy the commutation
relations
\be x_k x_l (P^-)^{kl}_{ij} =0, \eq
or equivalently
\be x_k x_l \R^{kl}_{ij}=q x_i x_j -\lambda \alpha (x \cdot x) g_{ij}, \eq
where $x \cdot x =x_k x_l g^{kl}=x_k x^k$ and
\be \label{alpha} \alpha=\frac{1}{1+q^{N-2}}. \eq
As a consequence the length
\be L=\alpha x \cdot x \eq
commutes with all the coordinates, $L x_i =x_i L$.

A calculus on quantum Euclidean space can be obtained by introducing
derivatives $\pa^i$ which
satisfy
\be \label{pxE} \pa^i x_j = \delta^i_j +q \R^{ik}_{jl} x_k \pa^l \eq
and
\be (P^-)^{ij}_{kl} \pa^l \pa^k =0. \eq
The Laplacian
\be \Delta=\alpha g_{ij} \pa^j \pa^i \eq
commutes with all derivatives, $\Delta \pa^i =\pa^i \Delta.$ One can define a
rescaling
operator
\be \label{Lambda} \Lambda =1+q \lambda x_i \pa^i +q^N \lambda^2 L \Delta, \eq
which satisfies
\be \label{Lambdax} \Lambda x_i =q^2 x_i \Lambda, \quad \pa ^i \Lambda =q^2
\Lambda \pa^i. \eq
A useful relation is
\be \label{pL} \pa^i L= q^2 L \pa^i +q^{2-N} x^i. \eq

The action of the vector fields $Z$ on $SO_q(N)$ induces in the standard way
an action on Euclidean space analogous to (\ref{Yx})
\be \label{Zx} Z^i_j x_k = x_m \R^{im}_{ln} Z^n_r \R^{lr}_{jk}. \eq
For $q$ real, the quantum Euclidean space can be endowed
with a reality structure as follows. For the coordinates one imposes the
reality
condition
\be \label{xbar} \overline{x_i}=g^{ij}x_j=x^i. \eq
Let us now define derivatives $\ph_i$ in terms of the conjugate derivatives by
\be \label{pconj} \ph_i=g_{ij} \ph^j=-q^N \overline{\pa^i}. \eq
The complex conjugate of Eq.(\ref{pxE}) can be transformed to the form
\be \label{phx} \hat{\pa}_i x^j= \delta^j_i +q^{-1} (\R^{-1})^{jl}_{ik} x^k
\hat{\pa}_l. \eq
The relation between the derivatives $\pa^i$ and their complex conjugates or
the $\ph_i$ can be
written\cite{OZ}
in the nonlinear form
\be \label{ph} \ph^i=\Lambda^{-1}(\delta^i_j +q^{N-1}\lambda \alpha x^i \pa_j)
\pa^j, \eq
which can be shown to satisfy Eq.(\ref{phx}).
Using Eq.(\ref{ph}), one can show that
\be \ph^i \pa^j =q \R^{ji}_{lk} \pa^k \ph^l. \eq

We wish to find a realization for the vector fields $Z$ of $SO_q(N,R)$ as
pseudodifferential operators
on Euclidean space. One way to find the appropriate expression is to proceed in
analogy with
Eq.(\ref{Uxp}) by writing similar terms but adjusting the coefficients so that
all relations
required of $Z$ are satisfied. It turns out that the correct formula is
\be \label{Zxp} Z^i_j=q^{-2} \delta^i_j +q^{-1}\lambda \pa^i x_j -q^{1-N}
\lambda x^i \ph_j -\lambda^2 L \pa^i  \ph_j. \eq
Using the relations given above for the calculus on Euclidean space, one can
verify that
$Z^i_j$ satisfies Eq.(\ref{Zx}) as well as
\be \label{pZ} \pa^i Z^j_k =\R^{ji}_{lm} Z^m_n \R^{ln}_{kr} \pa^r, \eq
and
\be \label{phZ} \ph^i Z^j_k =\R^{ji}_{lm} Z^m_n \R^{ln}_{kr} \ph^r. \eq

Combining Eqs.(\ref{Zx}), (\ref{pZ}) and (\ref{phZ}), one finds that $Z$
satisfies also
Eq.(\ref{ZZ}).
It is remarkable that $Z$, as given by Eq.(\ref{Zxp}) satisfies even the
orthogonality
relations Eqs.(\ref{gZZ}) and (\ref{ZZg}), without need for any further
normalization
as was necessary in Eqs.(\ref{ZY}) and (\ref{ZZU}).
This can be verified by direct computation and is due, apparently, to the fact
that
the $SO_q(N)$ $\R$ matrix already satisfies orthogonality relations.

Finally we may ask whether $Z$, as given by Eq.(\ref{Zxp}) satisfies the
natural
reality condition for $SO_q(N,R)$ which is
\be \label{Zreality} Z^\dagger =Z. \eq
It is very easy to see that this is indeed the case if one observes that
Eq.(\ref{Zxp})
can be written in the more symmetric form
\be q^2 Z^i_j= \delta^i_j +q \lambda \pa^i x_j +q \lambda \overline{x_i}
\overline{\pa^j} +\alpha q^N \lambda^2
\pa^i x_k \overline{x_k} \overline{\pa^j}, \eq
using Eq.(\ref{pL}).

On the other hand, if one does not impose Eq.(\ref{xbar})
and doesn't identify $\ph_i$, as given in Eq.(\ref{ph}), with the complex
conjugate derivative
$\overline{\pa^i}$ by Eq.(\ref{pconj}), then (\ref{Zreality}) will not be true.
However, Eq.(\ref{Zxp}) would still give a realization of
vector fields for the complex quantum group $SO_q(N)$ on Euclidean space.

In the differential calculus on a quantum space, one naturally introduces the
differentials of the
coordinates
\be \xi_i=dx_i. \eq
For quantum Euclidean space, they satisfy the commutation relations
\be \xi_k \xi_l (P^+)^{kl}_{ij}=0,\quad \xi_k \xi_l (P^0)^{kl}_{ij}=0, \eq
\be \label{xxi} x_i \xi_j = q \xi_k x_l \R^{kl}_{ij}, \eq
\be \pa^i \xi_j = q^{-1} (\R^{-1})^{ik}_{jl} \xi_k \pa^l. \eq
According to Eq.(\ref{xbar}) it is natural to introduce variables $\hat{\xi}_i$
related to $\overline{\xi_i}$ by
\be \overline{\xi_i}=g^{ij} \xih_j = \xih^i. \eq
The complex conjugate of Eq.(\ref{xxi}) can be written as
\be \label{xihx} \xih_i x_j =qx_k \xih_l \R^{kl}_{ij}. \eq
It was shown\cite{OZ} that the $\xih_i$ can be related to $\xi_i$ by a
(nonlinear) transformation
which was given explicitly there. It turns out that that transformation can be
written
very compactly as
\be \xih_i=\sigma q^N \Lambda \xi_k Z^k_i, \eq
where $\Lambda$ is given by Eq.(\ref{Lambda}). In this form one can easily
verify
that $\xih$ satisfies all desired relations. For instance Eq.(\ref{xihx})
follows immediately
from Eqs.(\ref{Lambdax}), (\ref{Zx}) and (\ref{xxi}). The requirement that
complex conjugation be an involution
restricts $\sigma$ to be a phase, $|\sigma|=1$. Vice versa, if one knows the
correct expression for $\xih_i$,
one can infer from it the formula for $Z^k_i$.

\section{Conclusion}
All above equations are "covariant". This means thay they go into themselves by
coaction transformations.
For instance, for all equations for $GL_q(N)$ from Eq.(\ref{RAA}) to
(\ref{Yp}), it is easy to see that
the transformation
\be A \rightarrow A B,\quad x \rightarrow x B, \quad \pa \rightarrow B^{-1}
\pa, \eq
\be Y \rightarrow B^{-1}Y B, \quad x_i \pa^i \rightarrow x_i \pa^i \eq
leaves them invariant. Here the matrix elements of $B$ are taken to commute
with
everything (which is the reason for using the word \underline{co}action) but
$B$ is itself a quantum matrix,
satisfying the analogue of Eq.(\ref{RAA}). It holds similarly for the complex
conjugate sector of $GL_q(N)$,
\be A^\dagger \rightarrow B^\dagger A^\dagger,\quad \xh \rightarrow B^\dagger
\xh, \quad \ph \rightarrow \ph(B^\dagger)^{-1}, \eq
\be Y^\dagger \rightarrow B^\dagger Y^\dagger (B^\dagger)^{-1}, \quad \ph_i
\xh^i \rightarrow \ph_i \xh^i \eq
(the relation $(B^\dagger)^{-1} =(B^{-1})^\dagger$ is used).
Analogous transformation laws leave invariant the
$SL_q(N), SO_q(N)$ equations as well as their respective real forms.

The realization of vector fields for $GL_q(N)$ and $SL_q(N)$ given in section 2
is equivalent to that given
earlier\cite{CSZ}. The formulas given here are simpler because of  a more
convenient choice of
notations and definitions. For instances, we use a right coaction and a
corresponding more convenient
lower index for the coordinates $x_i$ and upper index for the derivatives
$\pa^i$. The same applies
to a comparsion between the formulas written above for $SO_q(N)$ and earlier
ones\cite{OZ}. The reader
should have no difficulty in establishing the correspondence between the
conventions of these different
references.

A realization of vector fields for the orthogonal group in terms of
pseudodifferential operators on quantum Euclidean space has
been given by Gaetano Fiore\cite{Fiore}. He uses the explicit description of
the quantum Lie algebra
by Drinfeld and Jimbo, instead of Eqs.(\ref{ZZ}), (\ref{gZZ}) and (\ref{ZZg})
and gives explicit realizations for the vector fields
in that basis. Ours is an alternative solution of the same problem which has
perhaps the advantage of being
more symmetric and also covariant, as explained above.

\section{Acknowledgements}
This paper is our contribution to the Proceedings of the {\rm XX} International
Conference on
Group Theory Methods in Physics, Toyonaka, Japan (1994).

We wish to thank O. Ogievestky and M.Schlieker for many useful discussions and
comments.
This work was supported in part by the Director, Office of
Energy Research, Office of High Energy and Nuclear Physics, Division of
High Energy Physics of the U.S. Department of Energy under Contract
DE-AC03-76SF00098 and in part by the National Science Foundation under
grant PHY-90-21139.

\baselineskip 22pt



\begin{thebibliography}{99}

\bibitem{RTF}
N. Yu. Reshetikhin, L.A. Takhtajan and L.D. Faddeev,
{\it Quantization of Lie Groups and Lie Algebras},
Alg. i Anal. {\bf 1} 178 (1989) (Leningrad Math. J. {\bf 1} 193 (1990)).

\bibitem{Z1}
B. Zumino,
 {\it Introduction to the Differential Geometry of Quantum Groups},
 In: K. Schm\"udgen (ed.), Math. Phys. X,
 Proc. X-th IAMP Conf. Leipzig (1991),
 Springer-Verlag (1992).

\bibitem{Z2}
B. Zumino,
{\it Differential Calculus on Quantum Spaces
and Quantum Groups},
In: M.O., M.S., J.M.G. (eds.),
Proc. XIX-th ICGTMP Salamanca (1992),
CIEMAT/RSEF Madrid (1993).


\bibitem{SWZ1}
P. Schupp, P. Watts, and B. Zumino,
{\it Differential Geometry on Linear Quantum Groups},
Lett. Math. Phys. {\bf 25} 139 (1992).


\bibitem{SWZ2}
P. Schupp, P. Watts, and B. Zumino,
{\it Bicovariant Quantum Algebras and Quantum Lie Algebras},
Commun. Math. Phys. {\bf 157} 305 (1993).

\bibitem{WZ}
J. Wess and B. Zumino,
{\it Covariant Differential Calculus on the Quantum Hyperplane},
Nucl. Phys. B (Proc. Suppl.) {\bf 18B} 302 (1990).

\bibitem{GZB}
M.D.Gould, R.B. Zhang and A.J. Bracken,
{\it Generalized Gel'fand Invariants and Characteristic Identities
for Quantum Groups},
J. Math. Phys. {\bf 32} 2298 (1991).


\bibitem{Sch}
K. Schm\H{u}dgen,
{\it Operator Representation of the Real Twisted Commutation Relations},
J. Math. Phys. {\bf 33} 3211 (1994).

\bibitem{DSWZ}
B. Drabant, M. Schlieker, W. Weich and  B. Zumino,
{\it Complex Quantum Groups and their Quantum Enveloping Algebras},
Commun. Math. Phys. {\bf 147} 625 (1992).

\bibitem{CDSWZ}
C. Chryssomalakos, B. Drabant, M. Schlieker, W. Weich and B. Zumino,
{\it Vector Fields on Complex Quantum Groups},
Commun. Math. Phys. {\bf 147} 635 (1992).

\bibitem{Pod}
P. Podl\'{e}s,
{\it Complex Quantum Groups and their Real Representations},
Publ. R.I.M.S., Kyoto University {\bf 28} 709 (1992).

\bibitem{OZ}
O.Ogievetsky and B. Zumino,
{\it Reality in the Differential Calculus on q-Euclidean Spaces},
Lett. Math. Phys. {\bf 25} 121 (1992).

\bibitem{CSZ}
C. Chryssomalakos, P. Schupp and B. Zumino,
{\it Induced Extended Calculus on the Quantum Plane \/},
Alg. i Anal. {\bf 6} 252 (1994).

\bibitem{Fiore}
G.Fiore,
{\it Realization of $U_q(so(N))$ within the Differential Algebra on
${\bf R}_q^N$ }, preprint SISSA 90/93/EP, hep-th/9403033.

\end{thebibliography}
\end{document}